\documentclass[11pt,a4paper,twoside,english,notitlepage]{tran-l}
\usepackage{a4}
\usepackage{amsmath, amssymb}              
\usepackage{graphicx}
\usepackage{amsfonts}
\usepackage{babel}

\newtheorem{theorem}{Theorem}

\newtheorem{corollary}[theorem]{Corollary}

\newtheorem{definition}[theorem]{Definition}

\newtheorem{lemma}[theorem]{Lemma}

\newtheorem{remark}[theorem]{Remark}

\numberwithin{theorem}{section}

\def\0{{\bf 0}}

\def\R{{\bf R}}

\begin{document}

\title{Saari's Conjecture is True for Generic Vector Fields}

\author{Tanya Schmah}
\address{Department of Mathematics,
Division of ICS,
Macquarie University,
NSW 2109, Australia}
\email{schmah@maths.mq.edu.au}
\thanks{T. Schmah thanks the University of Surrey, the Bernoulli Centre at the Ecole Polytechnique F\'ed\'erale de 
Lausanne, and Wilfrid Laurier University, for their hospitality.
}

\author{Cristina Stoica}
\address{
Department of Mathematics,
Imperial College London,
SW7 2AZ London, UK;
and
Wilfrid Laurier University,
75 University Avenue W.,
Waterloo ON, N2L 3C5, Canada}
\email{cstoica@wlu.ca}
\thanks{
C. Stoica was supported by the MASIE (Mechanics and Symmetry in Europe)
Research Training Network of the European Union
(HPRN-CT-2000-0013) as a post-doctoral fellow at the University of Surrey,
and thanks Macquarie University for their hospitality.
}


\date{September 30, 2005}


\begin{abstract}
The simplest non-collision solutions of the $N$-body problem are the ``relative equilibria'', in which each body follows a circular orbit around the centre of mass and the shape formed by the $N$ bodies is constant. It is easy to see that the moment of inertia of such a solution is constant. In 1970, D. Saari conjectured that the converse is also true
for the planar Newtonian $N$-body problem:
relative equilibria are the only constant-inertia solutions. A computer-assisted proof for the 3-body case was recently given by R. Moeckel \cite{Moe05}.
 We present a different kind of answer: proofs that several generalisations of Saari's conjecture 
 are generically true.
Our main tool is jet transversality, including a new version
suitable for the study of generic potential functions.
\end{abstract}

\maketitle

\section{Introduction}
The Newtonian $N$-body problem concerns the motion of $N$ points
under the influence of a mutual gravitational force
\[
m_i \ddot{\mathbf{q}}_i 
= -\frac{\partial V}{\partial \mathbf{q}_i},
\]
for a potential function 
\[
V\left(\mathbf{q}_1,\dots,\mathbf{q}_N\right) 
= -\sum\limits_{1\le i<j\le N} \frac{m_i m_j}{\left\| \mathbf{q}_i - \mathbf{q}_j\right\|}
\]
(where $\mathbf{q}_i$ is position, $m_i$ is mass,
and the gravitational constant is taken to be $1$).
Without loss of generality, we will place the origin of the coordinate system at the centre of mass of the system.

Saari's conjecture concerns relative equilibria and the moment of inertia.
A \emph{relative equilibrium} of the $N$-body problem is a solution in which
the bodies move in circular orbits around the centre of mass, with
each body having the same constant angular velocity, so that the
shape formed by the $N$ points is constant.
For $N\ge 3,$ relative equilibria are the only explicit known periodic solutions
to the Newtonian $N$ body problem.
For $N=3$ it is known that there are only two kinds of relative equilibria: 
collinear (Eulerian) and 
equilateral triangle (Lagrangian).

The \emph{moment of inertia} of the $N$-body system is
\[
I=\sum\limits_{i=1}^N m_i \left\| \mathbf{q}_i \right\|^2
\]
(a factor of $\frac{1}{2}$ is often inserted).
It is a natural measure of the size of the system.
Physically, moment of inertia is a rotational analogue of mass:
note the similarity of $I$ to the kinetic energy 
$\displaystyle \frac{1}{2}\sum\limits_{i=1}^N m_i \left\| \dot{\mathbf{q}}_i \right\|^2 .$
Moment of inertia also has some interesting properties specific the Newtonian $N$-body problem,
arising from the homogeneity of $V$ and $I$. In particular,
if moment of inertia is conserved along any solution, then 
the potential and kinetic energies are separately conserved along
that solution.
For these reasons, moment of inertia is an interesting quantity to study in this problem.

Relative equilibria always have constant moment of inertia, since
each $\left\| \mathbf{q}_i \right\|$ remains constant.
In 1970 D. Saari conjectured that the converse is true:
\begin{verse}
\textbf{Saari's Conjecture:} \cite{Saa70} Every solution of the planar Newtonian $N$-body problem along which moment of inertia is conserved is a relative equilibrium.
\end{verse}
Attempts to prove this, by Saari and later by J. Palmore (1981), were unsuccessful.
Recent interest in Saari's conjecture is partly the indirect result of
the discovery of the``figure 8'' periodic solution to the three body problem,
numerically by C. Moore (1993), and analytically by A. Chenciner and R. Montgomery
\cite{CM00}.
Numerical calculations by C. Sim\'o indicated that this solution had nearly-constant
(but not constant) moment of inertia.

The conjecture has been proven for the planar $3$-body problem with equal masses by
C. McCord \cite{McC04} and J. Llibre and E. Pi\~na \cite{LlP02}.
F. Diacu,  E. P\'erez-Chavela, E.  and M. Santoprete \cite{DPS05}
have proven the conjecture for the collinear $N$-body problem.
Recently, R. Moeckel \cite{Moe05} has given a computer-aided proof for the 
general planar 3-body conjecture.
While this is very significant, interest in a simpler ``conceptual'' proof remains high.
An entry by A. Chenciner,
on a list of open problems compiled by K. C. Chen \cite{Ch03}, asks:
``Is there a conceptual proof for Saari's conjecture? Why not fix the moment of inertia tensor and ask the same question (maybe in higher dimensions)?''
The general case remains open.
Several researchers have worked to generalise Saari's conjecture
appropriately,
and to find counter-examples to these generalisations, including
G. Roberts \cite{Rob06}, M. Santoprete \cite{San04}, and
A. Hern\'andez-Gardu\~no, J.~K. Lawson, and J.~E. Marsden  \cite{HLM05}.

We have taken a different approach to the conjecture,
by asking whether
Saari's conjecture (now proven for N=3) is \emph{surprising}.
This question can be interpreted in the context of differential topology as:
is some appropriately generalised conjecture generically true?
We give positive answers to several such questions
using transversality theory.
Genericity and transversality theory are reviewed in Section \ref{secttrans},
which also contains a new version of jet transversality
suitable for the study of generic potential functions.

We generalise Saari's conjecture as follows.
Instead of the vector field given by the Newtonian potential, we consider arbitrary smooth
$G$-invariant vector fields $X$, for some Lie group $G$; 
and instead of moment of inertia we consider arbitrary smooth
$G$-invariant real-valued functions $F.$ 
We generalise Saari's conjecture to:
``the only solutions to $X$ along which $F$ is conserved are the relative equilibria''.
We show in Section \ref{sectgenF} that, for any given $X,$
the generalised conjecture is true for generic $F$
(Theorem \ref{genFsym}).
In Section \ref{sectgenX} we reverse the perspective and show that, for any given $F$
without too many critical points, the generalised conjecture is true for generic $X$
(Theorems \ref{gen1sym}).
We then prove analogous results
for generic Hamiltonian vector fields on symplectic
manifolds (Theorem \ref{Hamgen}), and for Hamiltonian vector fields for Hamiltonians of the form ``kinetic plus potential''
for a fixed kinetic energy and generic potentials (Theorem \ref{simpgen}).
While Theorems \ref{genFsym} and \ref{gen1sym} apply to general symmetry groups $G$,
our results for Hamiltonian vector fields assume that $G$ is trivial;
the corresponding equivariant conjectures remain open.

The basic idea in all of our proofs is the following: 
suppose $\gamma(t)$ is a solution to a vector field $X,$ 
and $F\circ \gamma$ is constant.
Then all derivatives of $F\circ \gamma$ are zero.
Equivalently, the Lie derivatives of $F$ with respect to $X,$ of all orders, are zero.
If $X$ is fixed, this puts restrictions on the  
Taylor expansion of $F,$ and we will show that generic functions $F$ do not
satisfy these restrictions, at any point in the phase space.
On the other hand, if $F$ is fixed, then the requirement that the derivatives be zero 
puts restrictions on the  
Taylor expansion of $X,$ and we will show that generic vector fields $X$ do not
satisfy these restrictions, at any point in the phase space.
The restrictions on the derivatives of $F$, or $X$, form a submanifold of an appropriate
jet space, to which we then apply jet transversality.

\section{Transversality} \label{secttrans}

In this section,
we summarise the tools from differential topology required to prove the
genericity results in the following sections.
All of the theorems stated here concern transversality. 
In order to state these results, we very briefly review genericity and jets.
The first half of this section is standard material,
a good reference being Hirsch \cite{H76}.
The main new result in the present section is Theorem \ref{newtrans},
which is a variation on jet transversality.
The proof uses Lemma \ref{globlem}, a globalisation lemma that we also use to
give a short proof of jet transversality for vector fields,
Theorem \ref{jettransvf}.

We consider only
$C^\infty$ paracompact, 
finite-dimensional manifolds without boundary.
Note that all such manifolds are second countable.
\emph{Smooth} will mean $\mathcal{C}^\infty.$
The function spaces $\mathcal{C}^\infty \left( M,N\right)$
will be given either the strong (``Whitney'') 
or the weak (``compact-open $\mathcal{C}^\infty$'') topology.
The choice of topology will be indicated either by subscripts ``$s$'' or ``$w$'',
or by the words ``strongly'' and ``weakly''.
Our main interest is in the strong topology; 
the weak topology will only be used in the proofs of Theorems
\ref{jettransvf} and \ref{newtrans} and in Lemma \ref{globlem}.

A subset is \emph{residual} if it contains the
intersection of a countable number of open dense subsets.
The spaces $\mathcal{C}^\infty \left(M,N\right)$,
with either the strong or the weak topology,
are $\emph{Baire spaces},$ 
meaning that  
all residual subsets are dense. 
A property is \emph{generic} within a given class of vector fields or functions, if
those with the given property form 
a residual subset of the class.
Loosely speaking, 
when a certain generic property is understood, a ``generic'' object means
one with that property.



\begin{definition}
Let $M$ and $N$ be manifolds, and $S$ a submanifold of $N.$ A map $f:M\to N$ is 
\emph{transverse} to $S$,
written $f\pitchfork S$,
 if, for every $x\in M$ such that $f(x)\in S,$ we have 
$T f \left(T_xM\right)+ T_{f(x)}S = T_{f(x)}N.$
\end{definition}


\begin{theorem}[Preimage]\label{preimage}
Let $f:M\to N$ be smooth and let $S$ be a codimension-$k$
submanifold of $N.$
If $f$ is transverse to $S$ then $f^{-1}(S)$ is either empty or a smooth codimension-$k$ 
submanifold of $M.$
\end{theorem}



\begin{theorem}[Elementary transversality]\label{elemtrans}
Let $M$ and $N$ be manifolds and let $S$ be a submanifold of $N.$
Then the smooth functions
$f:M\to N$ that are transverse to $S$ form a residual subset of $C_s^\infty\left(M,N\right).$ 
If $S$ is closed, then this subset is open dense
in $C_s^\infty\left(M,N\right)$. 
\end{theorem}

In the special case where $S$ has codimension larger than the dimension of $M,$
the previous two theorems combined imply that, for generic $f$,
the preimage $f^{-1}(S)$ is empty. 

In our application, we will need transversality to 
a set of restrictions on the Taylor series of a function or a vector field.
For this, we will need 
jet transversality (Theorems  \ref{jettrans} and \ref{jettransvf})
and a new related result (Theorem \ref{newtrans}).
The proofs of both of these use the following theorem,
which concerns genericity 
within a class of maps parametrised by a manifold $\mathcal{A}.$

\begin{theorem}[Parametric transversality]\label{partrans}
Let $M, N$ and
$\mathcal{A}$ be manifolds,
and let $S$ be a submanifold of $N$.
Let $\rho:\mathcal{A}\to \mathcal{C}_s^\infty\left(M,N\right)$ be a function
(not necessarily continuous).
Writing $\rho_\alpha = \rho\left(\alpha\right),$
define the evaluation map $\mathrm{ev}_\rho:\mathcal{A}\times M \to N$ by 
$\mathrm{ev}_\rho\left(\alpha,z\right) = \rho_\alpha(z).$
Suppose that $\mathrm{ev}_\rho$ is smooth and transverse to $S.$ 
Then the elements $\alpha$ of $\mathcal{A}$ such that 
$\rho_\alpha$ is transverse to $S$ 
form a residual subset of $\mathcal{A}.$
If $S$ is closed and $\rho$ is continuous
then the set of such $\alpha$ values is open dense.
\end{theorem}

We note that there is a version of this theorem due to R. Abraham \cite{A63,AR67}
that applies to Banach manifolds.
\medskip

We now briefly review \emph{jets}.
Let $M$ and $N$ be manifolds.
Two functions from $M$ to $N$ are \emph{equivalent at $x\in M$ up to order $k$}
if they have the same $k^{\textrm{th}}$ order Taylor expansion at $x$ in some coordinate charts. (This definition is coordinate-free.)
Such an equivalence class is called a \emph{$k$-jet} with \emph{source} $x$. 
The set of all such $k$-jets 
is written $J^k_x(M,N).$ The set
$J^k(M,N)$ is the union of these sets, for all $x.$
It is a smooth vector bundle over $M\times N,$ called the \emph{$k$-jet bundle}.
The $k$-jet of $f$ with source $x$ 
is written $j^kf(x).$ 
In local coordinates, $j^kf(x)$ ``is'' the $k^{\textrm{th}}$ order Taylor expansion of $f$ at $z.$ 
The \emph{$k$-jet extension} of $f:M\to N$ is the map 
\[
j^kf:M\longrightarrow J^k(M,N)\, ; \qquad
x\longmapsto j^kf(x)\, .
\]
If $f$ is smooth, then $j^kf$ is as well, so
there is a map
\begin{align} \label{E:jetext}
j^k:\mathcal{C}^\infty(M,N)\longrightarrow \mathcal{C}^\infty\left(M,J^k(M,N)\right),
\end{align}
taking $f$ to $j^kf$.
This map is continuous, with respect to the strong topologies on domain and codomain
\cite{GG73}.


\begin{theorem}[Jet transversality]\label{jettrans}
Let $M$ and $N$ be 
manifolds and let $S$ be a submanifold of $J^k(M,N).$ 
Then the set of functions
$f:M\to N$ such that $j^kf$ is transverse to $S$ is residual in 
$C^\infty_s\left(M,N\right)$, and open dense if $S$ is closed.
\end{theorem}

To apply jet transversality to vector fields, we need the
modified version in Theorem \ref{jettransvf}.
This result is known, but we are unaware of a proof in the literature.
We will prove it from Theorem \ref{jettrans}, 
using the globalisation technique in Lemma \ref{globlem}.
We will re-use the same globalisation
lemma in the proof of the new result in Theorem \ref{newtrans}.

Let $\mathfrak{X}^\infty(M)$ 
be the class of smooth vector fields on a manifold $M$.
Since $\mathfrak{X}^\infty(M)$ is a subset of $\mathcal{C}^\infty(M,TM)$,
it inherits a strong and a weak topology: the relative topology in each case.
It is easily verified that $\mathfrak{X}^\infty(M)$ is weakly closed
in $\mathcal{C}^\infty(M,TM)$. In fact, it is also weakly closed in
$\mathcal{C}^0(M,TM)$.
It follows that $\mathfrak{X}_s^\infty(M)$ (with the strong topology)
is a Baire space (see \cite{H76}).
Let $J^k\left(\mathfrak{X}^\infty(M)\right)$ be the subbundle
of $J^k\left(M,TM\right)$ consisting of all $k$-jets of vector fields.
The map
\[
j^k:\mathfrak{X}^\infty(M)\longrightarrow \mathcal{C}^\infty\left(M,J^k(\mathfrak{X}^\infty(M))\right)
\]
taking $X$ to $j^kX$,
is a restriction of the $k$-jet extension map in Equation (\ref{E:jetext}),
with $N=TM$, and hence the map
is continuous, with respect to the strong topologies on domain and codomain.

\begin{theorem}[Jet transversality for vector fields]\label{jettransvf}
Let $M$ be a 
manifold and let $S$ be a submanifold of 
$J^k\left(\mathfrak{X}^\infty(M)\right).$ 
Then the set of vector fields
$X\in \mathfrak{X}^\infty(M)$ such that $j^kX$ is transverse to $S$ is residual
in $\mathfrak{X}_s^\infty\left(M\right)$,
and open dense if $S$ is closed.
\end{theorem}

If $M$ is parallelisable, i.e. $TM$ is trivial, 
then this result follows directly from Theorem \ref{jettrans}.
For a general proof, we need a globalisation argument.
We now state such a result, Lemma \ref{globlem}, based on a similar one used by
Hirsch in the proof of jet transversality \cite{H76}.
Our result is for general vector bundles; the generality will allow us to 
re-use the argument in the proof of Theorem \ref{newtrans}.

For any smooth vector bundle $\pi:E\to M$, let $\Gamma^\infty(E)$ be the set
of smooth sections of $M$. Note that when $E=TM$, we have 
$\Gamma^\infty(E)= \mathfrak{X}^\infty(M)$.
Since $\Gamma^\infty(E)$ is a subset of $\mathcal{C}^\infty(M,E)$,
it inherits a strong and a weak topology: the relative topology in each case.
If $E$ is a trivial bundle,
then $E$ is isomorphic to $M\times \R^n$, for some $n$,
and there is a natural bijection between $\Gamma^\infty(M)$ and
$\mathcal{C}^\infty (M,\R^n)$.
The latter bijection is a homeomorphism with respect to either the strong topologies
on both spaces or the weak topologies on both spaces.


\begin{definition}
Let $\pi:E\to M$ be a smooth vector bundle.
A \emph{smooth class of sections} of $E$ is a family $\mathcal{X}$ of 
subsets $\mathcal{X}(U)\subseteq \Gamma^\infty(\pi^{-1}(U))$, defined for all 
open subsets $U\subseteq M$, 
satisfying the following ``localisation axioms'':
\begin{enumerate}
\item
If $\sigma \in \mathcal{X}(M)$ and
$U\subseteq M$ is open, then $\left.\sigma\right|_U \in \mathcal{X}(U)$.
\item
If $\sigma \in \Gamma^\infty\left(\pi^{-1}(U)\right)$ 
and there exists an open cover $\left\{U_i\right\}$ of $M$ 
such that $\left.\sigma\right|_{U_i}\in \mathcal{X}\left(U_i\right)$, for all $i$,
then $\sigma$ is in  $\mathcal{X}(M)$.
\end{enumerate}
\end{definition}

\begin{lemma}[Globalisation lemma] \label{globlem}
Let $\mathcal{X}$ be a smooth class of sections of a vector bundle $\pi:E\to M$.
If there exists an open cover $\mathcal{U}$ of $M$ such that,
for every open subset $U$ of an element of $\mathcal{U}$,
the set $\mathcal{X}(U)$ is strongly open
in $\Gamma^\infty(\pi^{-1}(U))$,
then $\mathcal{X}(M)$ is strongly open in $\Gamma^\infty(E)$.
If in addition, each $\mathcal{X}(U)$ is weakly dense in $\Gamma^\infty(\pi^{-1}(U))$,
then $\mathcal{X}(M)$ is also strongly dense in $\Gamma^\infty(E)$.
\end{lemma}

\begin{proof}
Cover $M$ with a countable 
and locally finite 
family of open sets $U_i$
such that the closure of each, $\overline{U_i}$, is
contained in an element $W_i$ of $\mathcal{U}$
(it suffices to let the $U_i$ be small enough coordinate discs, and then
apply paracompactness and second countability).
For each $i$, let 
\[
\mathcal{M}_i = 
\left\{ \omega \in \Gamma^\infty(E) : 
\left.\omega\right|_{U_i} \in \mathcal{X}\left(U_i\right)\right\}.
\]
From the localisation axioms, $\mathcal{X}(M)=\bigcap_i \mathcal{M}_i$.
By assumption, $\mathcal{X}\left(U_i\right)$ is strongly open
in $\Gamma^\infty\left(\pi^{-1}\left(U_i\right)\right)$.
This implies that each $\mathcal{M}_i$ is strongly open in $\Gamma^\infty(E)$.
Indeed, for any $\sigma \in \mathcal{M}_i$,
the restriction $\left.\sigma\right|_{U_i}$
has a strong basic neighbourhood $\mathcal{N}'_i$ in 
$\mathcal{X}\left(U_i\right)$, and the set
\[
\mathcal{N}_i = 
\left\{ \omega \in \Gamma^\infty(E) : 
\left.\omega\right|_{U_i} \in \mathcal{N}'_i \right\}
\]
is a strong basic neighbourhood of $\sigma$ in $\mathcal{M}_i$.
The neighbourhood $\mathcal{N}_i$ describes
restrictions on derivatives on some countable,
locally finite, family of compact subsets 
$L_{ij}\subseteq M$, all contained in $U_i$. 
Since the cover $\left\{U_i\right\}$ is countable and locally finite,
the family $\left\{L_{ij}\right\}$, for all $i$ and all $j$, is still countable and locally finite.
So $\mathcal{N}:=\bigcap_i \mathcal{N}_i$ is a strong basic neighbourhood 
of $\sigma$ in $\mathcal{X}(M)$.
Since $\sigma$ was arbitrary, this proves that $\mathcal{X}(M)$ is strongly open.

We now show that each $\mathcal{M}_i$ is strongly dense.
Let $\sigma \in \Gamma^\infty(E)$.
Let $\lambda_i:M\to \R$ be a $C^\infty$ function that equals $1$
on $\overline{U_i}$ and has compact support $K_i \subseteq W_i$.
For every $\tau\in \mathcal{X}\left(W_i\right)$,
define $F(\tau) \in \Gamma^\infty \left(E\right)$ by
\[
F(\tau) = \lambda_i \tau + \left(1-\lambda_i\right) \sigma
\]
(this is well-defined since $E$ is a vector bundle and $\lambda_i=0$
outside of the domain of $\tau$).
For every $\tau \in \mathcal{X}\left(W_i\right)$, the first localisation axiom
implies that $\left.\tau\right|_{U_i} \in \mathcal{X}\left(U_i\right)$, and then since
$\left.F(\tau)\right|_{U_i}=\left.\tau\right|_{U_i}$, we have $F(\tau)\in \mathcal{M}_i$.
Since $\lambda_i$ has compact support, all of its derivatives are bounded, 
so $F$ is a weakly continuous map,
\[
F:\left( \mathcal{X}\left(W_i\right) \subseteq \Gamma^\infty_w \left(\pi^{-1}(W_i)\right) \right) \longrightarrow
\mathcal{M}_i\subseteq \Gamma^\infty_w \left(E\right)\, .
\]
Now every neighbourhood of $\sigma$ in  $\Gamma_s^\infty(E)$
(with the strong topology)
contains a strong basic neighbourhood
of the form $\mathcal{N}=\bigcap_j\mathcal{N}_j$,
where each $\mathcal{N}_j$ restricts derivatives on some 
set $L_j\subseteq M$,
and the family $\left\{L_j\right\}$ is a locally finite.
By the local finiteness, the compact set $K_i$
has nontrivial intersection with
only a finite number of the sets $\mathcal{L}_j$.
Let $\mathcal{N}_w$ be the intersection of the corresponding sets $\mathcal{N}_j$,
and note that since the intersection is finite, $\mathcal{N}_w$ is open in the
weak topology. 
For every $\tau\in \Gamma^\infty \left(\pi^{-1}(W_i)\right)$,
since $F(\tau) - \sigma$ has compact support $K_i$,
it follows that $F(\tau)\in \mathcal{N}$
if and only if  $F(\tau)\in \mathcal{N}_w$.
Thus, to prove that $\sigma$ is in the strong closure of $\mathcal{M}_i$,
it suffices to find local sections 
$\tau \in \mathcal{X}\left(W_i\right)$
such that $F(\tau)$ is arbitrarily close to $\sigma$ in the weak topology.
Since $F$ is weakly continuous, it suffices
to find $\tau \in \mathcal{X}\left(W_i\right)$ arbitrarily close to 
$\left.\sigma\right|_{U_i}$ in the weak topology.
But we have assumed that $\mathcal{X}\left(W_i\right)$ is weakly dense
in $\Gamma^\infty\left(\pi^{-1}\left(W_i\right)\right)$,
which finishes the proof that $\sigma$ is in the strong closure of $\mathcal{M}_i$.
Hence $\mathcal{M}_i$ is strongly dense.

We have shown that each $\mathcal{M}_i$ is strongly open and strongly dense 
in $\Gamma_s^\infty(E)$.
Recall that  $\mathcal{X}(M)=\bigcap_i \mathcal{M}_i$.
Since this intersection is countable and $\Gamma_s^\infty(E)$ is a Baire space,
it follows that $\mathcal{X}(M)$ is strongly dense. 
\end{proof}



\begin{proof}[Proof of Theorem \ref{jettransvf}
 (``Jet transversality for vector fields'')]
Let $M$ be a manifold and $S$ a submanifold of $J^k\left(\mathfrak{X}^\infty(M)\right)$.
First, suppose that $S$ is closed.
For every open $U\subseteq M$, 
note that $J^k\left(\mathfrak{X}^\infty(U)\right)$) is an open subset of  
$J^k\left(\mathfrak{X}^\infty(M)\right)$).
Let $\pitchfork(U)$ be the the set of all 
$X\in \mathfrak{X}^\infty(U)$ such that $j^kX$ is transverse to 
$S\cap J^k\left(\mathfrak{X}^\infty(U)\right)$.
It is easily verified that the family of sets $\pitchfork(U)$, for all open $U\subseteq M$
is a smooth class of sections of $TM$.
If $TU$ is a trivial bundle, 
then $\mathcal{X}_s^\infty(U)$ is homeomorphic
to $\mathcal{C}_s^\infty\left(U,\R^n\right)$, where $n$ is the dimension of $M$,
and so it follows from Theorem \ref{jettrans} that $\pitchfork(U)$ is 
open dense in $\mathcal{X}_s^\infty(U)$.
Since strongly dense implies weakly dense,
and all bundles are locally trivial,
the class $\left\{\pitchfork(U)\right\}$ satisfies the conditions of 
the globalisation lemma (Lemma \ref{globlem}), 
with $\mathcal{U}$ being the set of $U\subseteq M$ such that $TU$ is trivial.
Hence
$\pitchfork(M)$ is strongly dense in $\mathfrak{X}^\infty(M)$.

The set $\pitchfork(M)$ is the preimage by $j^k$
of
\[
\left\{
f\in \mathcal{C}^\infty\left(M, J^k\left(\mathfrak{X}^\infty(M)\right)\right)
: f \pitchfork S
\right\},
\]
which by Theorem \ref{jettrans} is strongly open.
Since $j^k$ is strongly continuous,
it follows that $\pitchfork{M}$ is strongly open. 

If $S$ is not closed, express it as the countable union of 
closed coordinate disks $S_i$. 
For every $i$, let $\pitchfork_i(M)$ be the set of all smooth vector fields
such that $j^kX$ is transverse to $S_i$.
Since $\pitchfork(M) = \bigcap_i \pitchfork_i(M)$,
and we have shown that each of the sets $\pitchfork_i(M)$ is open dense
(in the strong topology),
it follows that $\pitchfork(M)$ is residual.
\end{proof}

The jet transversality results in Theorems \ref{jettrans} and \ref{jettransvf}
 will be sufficient to prove Theorems \ref{genF}, \ref{gen1} and \ref{Hamgen}.
For Theorem \ref{simpgen}, concerning generic potentials, a new variation on jet transversality is required.
We consider
vector fields parametrised by $\mathcal{C}^\infty\left(Q,\R\right),$
for some manifold $Q$,
and a map $\rho$ from $\mathcal{C}^\infty\left(Q,\R\right)$ to 
$\mathcal{C}^\infty\left(P,J^k\left(\mathfrak{X}^\infty(P)\right)\right)$,
where $P$ is a vector bundle over $Q.$
We have in mind $P=T^*Q$ and
the elements of $\mathcal{C}^\infty\left(Q,\R\right)$
being potential functions.
We assume that $\rho(f)$ depends only on the first $m$
derivatives of $f$, i.e. that
$\left(\rho(f)\right)(z)$ depends only on $j^mf\left(\pi(z)\right)$ and $z$.
Since, for any $m,$ the set $J^m\left(Q,\R\right)$ is a bundle over $Q\times \R,$ and hence over $Q,$
we may form the
Whitney sum $J^m\left(Q,\R\right)\oplus P,$ which is a bundle over $Q.$

\begin{theorem}[``Jet transversality for potentials''] \label{newtrans}
Let $\pi:P\to Q$ be a 
smooth vector bundle, and let
$S$ be a submanifold of $J^k(\mathfrak{X}^\infty(P)).$ 
Let $
\rho:\mathcal{C}^\infty\left(Q,\R\right)\to
\mathcal{C}^\infty\left(P, J^k\left(\mathfrak{X}^\infty(P)\right)\right)
$
be a function.
Suppose that the map
\begin{align*}
\Phi: J^m(Q,\R) \oplus P &\longrightarrow  J^k\left(\mathfrak{X}^\infty(P)\right) \\
\left(j^mf\left(\pi(z)\right),z\right) &\longmapsto \left(\rho(f)\right)(z)
\end{align*}
is well-defined, smooth and transverse to $S$.
Then the set of $f\in \mathcal{C}^\infty\left(Q, \R\right)$
such that $\rho(f)$ is transverse to $S$ 
is residual in
$C^\infty_s\left(Q,\R\right)$,
and open dense if $S$ is closed.
\end{theorem}

Our proof is based on Hirsch's proof of jet transversality \cite{H76}. 
The main idea is: starting with any potential $f$, consider the set of all polynomial
perturbations of $f$, which is a finite-dimensional vector space, and apply 
parametric transversality to show that the subset of transverse perturbations is dense.
This idea suffices locally, and the global result then follows
from Lemma \ref{globlem} (globalisation), applied to
the trivial bundle $Q\times \R \to Q$
(sections of which ``are'' elements of $\mathcal{C}^\infty(Q,\R)$).

\begin{proof}[Proof of Theorem \ref{newtrans}
(``Jet transversality for potentials'')]
Let 
\[
\pitchfork(Q) = \left\{ f\in \mathcal{C}^\infty(Q,\R) :
\rho(f) \pitchfork S
\right\}.
\]
First, suppose that 
$Q$ is an open subset of $\R^n$, and 
$P=Q\times \R^{p-n}$ (with $p\ge n$), and that $S$ is closed.
Let $f\in \mathcal{C}^\infty(Q,\R)$. We will show that $f$ is in the weak closure of
$\pitchfork(Q)$.
Each element of $J_{\mathbf{0}}^m\left(Q,\R\right)$ can be expressed uniquely as
$j^mg(\mathbf{0})$ for an $m^{\mathrm{th}}$-order polynomial $g$ on $Q=\R^n$,
so the following map is well defined:
\begin{align*}
\alpha:J_{\mathbf{0}}^m\left(Q,\R\right) &\longrightarrow \mathcal{C}^\infty(Q,\R) \\
j^mg(\mathbf{0}) &\longmapsto f + g\, ,
\end{align*}
for every $m^{\mathrm{th}}$-order polynomial $g$.
This map is weakly continuous,
and $\alpha(\mathbf{0}) = f$,
so it suffices to show that $\mathbf{0}$ is in the weak closure of
$\left\{ x \in J_0^m\left(Q,\R\right) : \rho(\alpha (x)) \pitchfork S\right\}$.
We will accomplish this by applying parametric transversality (Theorem \ref{partrans}) to
the map $F:= \rho\circ \alpha$, which will prove that the set 
$\left\{ x \in J_0^m\left(Q,\R\right) : F(x) \pitchfork S\right\}$
is strongly dense, and hence weakly dense, in  $J_0^m\left(Q,\R\right)$.
To do this, we must show that $\mathrm{ev}_F$, the evaluation map of $F$, 
is smooth and transverse to $S$.

Let $\Phi$ be as in the statement of the theorem. Then, for every $z\in P$,
\[
\left(F(j^mg(\mathbf{0}))\right)(z) = \rho(f+g)(z) = \Phi\left(j^m (f+g)\left(\pi(z)\right), z\right) ,
\]
so we can factor $\mathrm{ev}_F$ as follows (this defines $G$):
\[
\begin{array}{cccccc}
\mathrm{ev}_F:&J_{\mathbf{0}}^m\left(Q,\R\right) \times P 
& \overset{G}{\longrightarrow} &J^m(Q,\R) \oplus P
&\overset{\Phi}{\longrightarrow} &J^k\left(\mathfrak{X}^\infty(P)\right) \\
&\left(j^mg(\mathbf{0}),z\right) & \longmapsto 
&\left(j^m(g+f)\left(\pi(z)\right),z\right) & \longmapsto &\rho(g+f)(z)
\end{array}
\]
By assumption, $\Phi$ is smooth and transverse to $S$,
so it suffices to show that $G$ is a smooth submersion.

Since we are assuming $Q$ is an open subset of $\R^n$,
there is a natural isomorphism
\begin{align*}
J^m(Q,\R) &\longrightarrow J^m_{\mathbf{0}}(Q,\R) \times Q \\
j^mf(q) & \longmapsto \left(j^m\left(f+\Sigma_q\right)(\mathbf{0}), q\right)\, ,
\end{align*}
where $\Sigma_q$ is the shift map defined by $\Sigma_q(q') = q'+q$.
Since we are also assuming $P=Q\times \R^{p-n}$, 
the above isomorphism induces one between the Whitney sum
$J^m(Q,\R) \oplus P$ and
$J^m_{\mathbf{0}}(Q,\R) \times P$. Making this identification, $G$ can be written as
\begin{align*}
G:J^m_{\mathbf{0}}(Q,\R) \times P & \longrightarrow J^m_{\mathbf{0}}(Q,\R) \times P \\
\left(j^mg(\mathbf{0}), z\right) 
&\longmapsto \left(j^m\left((g+f)\circ \Sigma_{\pi(z)}\right)(\mathbf{0}),z\right)\, .
\end{align*}
This is a smooth map, since $f$ and $g$ are smooth, for all polynomials $g$.
For any fixed $z$, the map
\[
j^mg(\mathbf{0}) 
\longmapsto 
j^m\left((f+g)\circ \Sigma_{\pi(z)}\right)(\mathbf{0})
= j^m\left(f\circ \Sigma_{\pi(z)}\right)(\mathbf{0})
+j^m\left(g\circ \Sigma_{\pi(z)}\right)(\mathbf{0})
\]
is an affine bijection on $J^m(Q,\R)$,
since $g\circ \Sigma_{\pi(z)}$ is a polynomial with coefficients that are 
linear functions of the coefficients of $g$. 
It follows that $G$ is a smooth submersion.
Thus $\mathrm{ev}_F$ is smooth and transverse to $S$.
As explained above, this implies that $f$ is in the weak closure of
$\pitchfork(Q)$.
Since this argument holds for any $f\in \mathcal{C}^\infty(Q,\R)$, we have shown
that $\pitchfork(Q)$
is weakly dense in $\mathcal{C}^\infty(Q,\R)$.

Still assuming $P=Q\times \R^{p-n}$, 
let $\pi_2:P\to \R^{p-n}$ be projection onto the second factor.
Define
\begin{align*}
H:\mathcal{C}_s^\infty\left(Q,J^m(Q,\R)\right)
&\longrightarrow \mathcal{C}_s^\infty\left(P,J^m(Q,\R) \times \R^{p-n}\right)\\
\varphi &\longmapsto (\varphi \circ \pi) \times \pi_2\, ,
\end{align*}
where the $\times$ operation is defined by
$\left((\varphi \circ \pi) \times \pi_2\right) (z) = \left(\varphi \circ \pi (z), \pi_2(z)\right)$.
Using well-known facts about the strong topology \cite{GG73},
it can be checked that $H$ is strongly continuous.
Identifying $J^m(Q,\R) \times \R^{p-n}$ with $J^m(Q,\R)\oplus P$,
the map $H$ is given by 
$H(\varphi)(z) = \left(\varphi \circ \pi(z), z\right) \in J^m(Q,\R)\oplus P$.
For every $f \in \mathcal{C}_s^\infty\left(Q,J^m(Q,\R)\right)$ and every $z\in P$,
\[
\left(\rho(f)\right)(z) = \Phi\left(j^mf \circ \pi(z),z\right)
= \left(\Phi \circ H \circ j^m(f)\right)(z) \, .
\]
Hence $\rho = \Phi \circ H \circ j^m $. 
Since $j^m$ is strongly continuous, it follows that
$\rho$ is strongly continuous.
Since
\[
\pitchfork(Q) = \rho^{-1} 
\left\{\varphi\in \mathcal{C}^\infty\left(P, J^k\left(\mathfrak{X}^\infty(P)\right)\right)
: \varphi \pitchfork S\right\},
\]
Theorem \ref{elemtrans} (elementary transversality) implies that
$\pitchfork(Q)$ is strongly open.

Now let $\pi:P\to Q$ be any smooth vector bundle. Again, we assume $S$ is closed. 
For every open $U\subseteq Q$, we define a ``localisation'' of $\rho$ by 
$\rho_U(h)(z) = \Phi\left(j^mh\left(\pi(z)\right),z\right)$,
\[
\rho_U:\mathcal{C}^\infty\left(U,\R\right)\longmapsto
\mathcal{C}^\infty\left(\pi^{-1}(U), J^k\left(\mathfrak{X}^\infty(\pi^{-1}(U))\right)\right),
\]
and define
\[
\pitchfork(U) = \left\{ h \in \mathcal{C}^\infty(U,\R) 
: \rho_U(h) \pitchfork \left(S\cap J^k\left(\mathfrak{X}^\infty(\pi^{-1}(U))\right)\right)\right\}.
\]
It is easily verifed that the family $\left\{ \pitchfork(U)\right\}$ is a smooth
class of sections of the trivial bundle $Q\times \R \to Q$.
If the restriction of $\pi:P\to Q$ to $U$ is trivial,
our earlier arguments show that $\pitchfork(U)$ is strongly open
and weakly dense in 
$\mathcal{C}^\infty\left(U,\R  \right)$.
Since all bundles are locally trivial,
$Q$ can be covered by such open sets $U$.
By Lemma \ref{globlem} (globalisation), $\pitchfork(Q)$ is
strongly open and strongly dense in $\mathcal{C}^\infty\left(Q,\R  \right)$.

If $S$ is not closed, express it as the countable union of 
closed coordinate disks $S_i$. 
For every $i$, let $\pitchfork_i(Q)$ be the set of all $f\in \mathcal{C}^\infty(Q,\R)$
such that $j^kf$ is transverse to $S_i$.
Since $\pitchfork(Q) = \bigcap_i \pitchfork_i(M)$,
and we have shown that each of the sets $\pitchfork_i(M)$ is open dense
(in the strong topology),
it follows that $\pitchfork(M)$ is residual.

\end{proof}

\begin{remark}[Transversality to families of submanifolds]
Since the countable intersection of residual sets is residual, any of the above theorems
may be applied repeatedly to each of a finite or countable family of submanifolds,
yielding a residual subset of maps transverse to all of the submanifolds.
\end{remark}

\section{Generic quantities are never conserved by a given vector field} \label{sectgenF}

Consider a smooth vector field $X$ on a manifold $P,$ and a smooth real-valued 
function $F:P\to \R$.
We will begin by showing that, for a given function $F,$
generic vector fields $X$ have no non-equilibrium solutions conserving $F$ (Theorem \ref{gen1}).
Theorem \ref{gen1sym} is the ``equivariant'' version of this result: if $F$ is $G$-invariant, for some
$G$ acting freely, then generic $G$-invariant vector fields  have no solutions conserving $F$
other than the relative equilibria.

Recall that a property is \emph{generic} within a given space of functions
if those functions with the property form a residual subset of the space;
in fact, in the following theorem, the subset is open dense. 
The strong function space topologies will be assumed in this section and
the remainder of this article.

\begin{theorem}\label{genF}
Let $P$ be a 
manifold and let
$X$ be a smooth vector field on $P.$
Then, for generic smooth functions
$\displaystyle{F:P\to\mathbf{R}},$
the only solutions to $X$ along which $F$ is conserved are the equilibria (i.e. fixed points).
\end{theorem}

\begin{proof}
Let $n$ be the dimension of $P,$ and
consider
a solution $\gamma(t)$ of $X.$
For any smooth function $F,$ the derivatives $\left(F\circ \gamma\right)^{(k)}(0),$
for $k=1,\dots,n+1,$ are determined by $X$ and by
$j^{n+1}F\left(\gamma(0)\right).$
Hence we can define a map
$\Psi:J^{n+1}(P,\R) \to \R^{n+1}$ by
\[
\Psi\left(j^{n+1}F(z)\right) =
\left(\left(F\circ \gamma\right)'(0),\dots,\left(F\circ \gamma\right)^{(n+1)}(0)\right),
\]
where $\gamma$ is the solution of $X$ with initial condition $\gamma(0)=z.$
If $F$ is constant along $\gamma$ then
all derivatives $\left(F\circ \gamma\right)^{(i)}(t)$ must be zero, for all $t,$ so
$\gamma(t)$ must remain in $\left(\Psi\circ j^{n+1}F\right)^{-1}(\0).$
We will show that, for generic $F,$ this set is empty except for the equilibrium points of $X.$

Let $S_e\subseteq J^{n+1}\left(P,\R\right)$
be the set of all $j^{n+1}F(z)$ for which $X(z)=0,$
i.e. $z$ is an equilibrium point of the flow of $X.$
Let $S_1$ be the complement of $S_e$ in $J^{n+1}\left(P,\R\right),$
and note that $S_1$ is open 
in $J^{n+1}\left(P,\R\right).$
We will show that, for generic $F,$
the set $\left(\left.\Psi\right|_{S_1}\right)^{-1}(\0)$
is either empty or a codimension-$(n+1)$ submanifold of
$S_1,$ and hence of $J^{n+1}\left(P,\R\right)$ as well.
To do this it suffices to show that $\left.\Psi\right|_{S_1}$ is a submersion.

We begin by computing $\Psi\left(j^{n+1}F(z)\right)$ in local coordinates for $P.$
Using superscripts to denote components of $X$ and subscripts to denote partial differentiation, using the summation convention for repeated indices,
and evaluating all derivatives at the same point $z=\gamma(0),$
we obtain:
\begin{align*}
\left(F\circ \gamma\right)'(0) &= F_{j} X^{j} \\
\left(F\circ \gamma\right)''(0) &= F_{j_1j_2} X^{j_1} X^{j_2} + F_{j_1} X^{j_1}_{j_2} X^{j_2} \\
&\vdots \\
\left(F\circ \gamma\right)^{k}(0) &=
F_{j_1j_2\cdots j_{k}} X^{j_1} X^{j_2}\cdots X^{j_{k}}
+ \textrm{(terms of lower order in $F$)}
\end{align*}
We write the components of $\Psi$ as $\Psi_k,$ for $k=1,\dots,n+1,$ so
$\Psi_k(j^{n+1}F(z)) = \left(F\circ \gamma\right)^{(k)}(0).$
We now fix a $j^{n+1}F\left(z\right)\in S_1$
and consider $D\Psi$ at this point.
The partial derivatives of $F$ up to order $n+1$ are coordinates for
$J^{n+1}(P,\R)$, the domain of $\Psi.$
We will choose $n+1$ of them and use them to show that 
the rank of 
$D\Psi\left(j^{n+1}F\left(z\right)\right)$ is $n+1.$
By definition of $S_1,$ we see that $z$ is not a critical point of $X,$
so there exists a $j\le n$ such that $X^j(z)\ne 0.$
It follows that, for every $k,$
the coefficient of $F_{jj\cdots j}$ in $G_k$ is nonzero, where
$F_{jj\cdots j}$ means the $k^{th}$-order partial derivative of $F$ with respect to
the $j^{th}$ component of $z.$
This shows that the partial derivative of $\Psi_k$ with respect to $F_{jj\cdots j}(z)$ is nonzero.
Since the $k^\mathrm{th}$ order derivative $F_{jj\cdots j}$ does not appear in $\Psi_1,\cdots, \Psi_{(k-1)},$ and this argument holds for all $k=1,\dots,n+1,$ this shows that
the rank of  $D\Psi\left(j^{n+1}F\left(z\right)\right)$ is $n+1.$
Since the
argument holds for all $j^{n+1}F(z)\in S_1,$ this shows that
the restriction of $\Psi$ to $S_1$ is a submersion.

It follows that $S_2:=\left(\left.\Psi\right|_{S_1}\right)^{-1}(\0)$ is a closed smooth
submanifold of $S_1,$
and hence of $J^{n+1}\left(P,\R\right)$ as well, and that $S_2$
is either empty or of codimension $(n+1).$
Note that $\Psi^{-1}(\mathbf{0})\subseteq S_e \cup S_2.$
It follows that
\[
\left(\Psi\circ j^{n+1}F \right)^{-1}(\0) \subseteq \left(j^{n+1}F\right)^{-1}(S_e \cup S_2).
\]

Let $\mathcal{B}$ be the set of all
$F\in \mathcal{C}^\infty(P,\R)$ such that
$j^{n+1}F$ is transverse to $S_2.$
By jet transversality (Theorem \ref{jettrans}),
$\mathcal{B}$ is an open dense subset of $\mathcal{C}^\infty(P,\R).$
For any $F\in \mathcal{B},$
since $S_2$ is either empty or has codimension $n+1,$
its preimage $\left(j^{n+1}F\right)^{-1}(S_2)\subseteq P$ must be empty.
As noted earlier, any nontrivial solution of $X$ along which $F$ is constant
must remain in $\left(\Psi\circ j^{n+1}F \right)^{-1}(\0).$ Since we have just
shown that $\left(j^{n+1}F\right)^{-1}(S_2)$ is empty, the only possible solutions are those that
remain in $\left(j^{n+1}F\right)^{-1}(S_e).$ But this is the set of equilibrium points of $X.$
Hence the only solutions of $X$ along which $F$ is constant are equilibria.
\end{proof}

Since the $N$-body problem has no equilibrium solutions, we have the following:

\begin{corollary}
In the Newtonian $N$-body problem (planar or spatial), for generic smooth real-valued functions $F$
on phase space, there are no 
solutions along which $F$ is conserved.
\end{corollary}

The fact that relative equilibria conserve moment of inertia 
implies that moment of inertia is not generic in the sense of the corollary. 
This non-genericity
is obviously related to the symmetry of the problem:
in the case of the planar problem in center-of-mass coordinates,
there is an $SO(2)$ symmetry on the phase space $T^*Q\cong\R^{4N-4}$
that conserves the vector field and the moment of inertia.
Thus when generalising Saari's conjecture,
it is most natural to restrict attention to $SO(2)$-symmetric functions $F,$ for which
the conjecture states:
``a solution of the planar $N$-body problem
conserves $F$ if and only if the solution is a relative equilibrium''.
We now ask: is this generalisation of Saari's conjecture true
for generic $SO(2)$-invariant functions $F$?
For simplicity we consider only free symmetries.
Of course, the $SO(2)$ symmetry in the planar $N$-body problem is not free,
but it is free if we remove the origin (the centre of mass)
from the configuration space, which
has no effect on Saari's conjecture
since this configuration is an $N$-body collision, at which point  the potential is
undefined.

If a Lie group $G$ acts freely, properly and smoothly on a manifold $P,$
then $P/G$ is a smooth manifold.
The class $\mathfrak{X}^G(P)$ of smooth $G$-invariant vector fields on $P$
is isomorphic to $\mathfrak{X}^\infty(P/G).$
Similarly, the class $\mathcal{C}^G(P,\R)$  of smooth 
$G$-invariant functions on $P$
is isomorphic to $\mathcal{C}^\infty(P/G,\R),$ with this isomorphism defining the
Whitney (strong) topology on  $\mathcal{C}^G(P,\R).$
Given an $X\in \mathfrak{X}^G(P),$ we apply
Theorem \ref{genF} to the corresponding vector field $\bar{X}\in \mathfrak{X}\left(P/G\right),$ 
concluding that for generic $\bar{F}\in \mathcal{C}^\infty(P/G,\R),$
there are no non-equilibrium solutions to $\bar{X}$ that conserve $\bar{F}.$
The corresponding conclusion in the original phase space is:

\begin{theorem}\label{genFsym}
Let $P$ be a 
manifold, let $G$ be a Lie group acting smoothly, properly and freely on $P,$ and let
$X$ be a smooth $G$-invariant vector field on $P.$
Then, for generic smooth $G$-invariant functions
$\displaystyle{F:P\to\mathbf{R}},$
the only solutions to $X$ along which $F$ is conserved are the relative equilibria.
\end{theorem}

\begin{corollary} In the planar Newtonian $N$-body problem,
for generic smooth $SO(2)$-invariant functions
$\displaystyle{F:P\to\mathbf{R}},$
the only solutions along which $F$ is conserved are the relative equilibria.
\end{corollary}


In the next section, we reverse the point of view by fixing a function $F$ and considering generic vector fields.

\section{Generic vector fields never conserve a given quantity} \label{sectgenX}

Consider a vector field $X$ on a manifold $P,$ and a smooth function $F:P\to \R.$
As earlier, we study the question of whether there is at least one solution to $X$ along which $F$
is conserved, but this time we fix $F$ and allow $X$ to vary.
Note that, if $F$ is constant, it is trivially true that all solutions conserve $F$. 
We will exclude this case, and more generally, require that the critical
points of $F$ be contained in some codimension-$1$ manifold.

We will begin by showing that, for a given function $F,$
generic vector fields have no non-equilibrium solutions conserving $F$ (Theorem \ref{gen1}).
We next prove the equivariant version of this result: if $F$ is $G$-invariant, for some
$G$ acting freely, then generic $G$-invariant vector fields  have no solutions conserving $F$
other than the relative equilibria (Theorem \ref{gen1sym}).
Theorem \ref{Hamgen} concerns Hamiltonian vector fields on symplectic manifolds, and
Theorem \ref{simpgen} concerns Hamiltonian vector fields for Hamiltonians
of the form ``kinetic plus potential'' for a fixed kinetic energy.


\begin{theorem}\label{gen1}
Let $P$ be a 
manifold. Let
$\displaystyle{F:P\to\mathbf{R}}$ be smooth and suppose that
its critical points
are contained in a codimension-$1$ submanifold of $P.$
Then
there exists a residual subset $\mathcal{C}$ of $\mathfrak{X}^\infty(P)$
such that no vector field in 
$\mathcal{C}$ has any non-equilibrium
solution along which $F$ is constant,
and the equilibrium solutions of any vector field in $\mathcal{C}$ are isolated.
\end{theorem}

\begin{remark} The following proof is very similar to that of Theorem \ref{genF}.
The main differences are due to that fact that the terms in $\Psi_k$ containing the
 highest-order partial derivatives of $X$ contain factors of $F_i$ as well as $X^i.$
\end{remark}

\begin{proof}
Let $n$ be the dimension of $P.$
Consider a vector field $X\in \mathfrak{X}^\infty\left(P\right)$
and a solution $\gamma(t)$ of $X.$
The derivatives $\left(F\circ \gamma\right)^{(i)}(0),$
for $i=1,\dots,n+1,$ are determined by 
$F$ and $j^{n}X\left(\gamma(0)\right).$ Hence
we can define a map
$\Psi:J^{n}(\mathfrak{X}^\infty(P)) \to \R^{n+1}$ by
\[
\Psi\left(j^{n}X(z)\right) =
\left(\left(F\circ \gamma\right)'(0),\dots,\left(F\circ \gamma\right)^{(n+1)}(0)\right),
\]
where $\gamma$ is the solution of $X$ with initial condition $\gamma(0)=z.$
If $F$ is constant along $\gamma$ then 
all derivatives $\left(F\circ \gamma\right)^{(k)}(t)$ must be zero, for all $t,$ so
$\gamma(t)$ must remain in $\left(\Psi\circ j^nX\right)^{-1}(\0).$ 

Let $S_e$ be the codimension-$n$ submanifold of $J^n\left(\mathfrak{X}^\infty(P)\right)$
consisting of all $j^{n}X(z)$ for $X$ such that $X(z)=0,$
i.e. $z$ is an equilibrium point of the flow of $X.$
By assumption, the critical points of $F$ are all contained in some single codimension-$1$ 
manifold $Z.$
Let $S_F$ be the codimension-$1$ submanifold of $J^n\left(\mathfrak{X}^\infty(P)\right)$
consisting of all $j^nX(z)$ such that $z\in Z.$
Let $S_1$ be the complement of $S_e\cup S_F$ in $J^n\left(\mathfrak{X}^\infty(P)\right),$
and note that $S_1$ is open dense in $J^n\left(\mathfrak{X}^\infty(P)\right).$
We will show that for, generic $X,$
the set $\left(\left.\Psi\right|_{S_1}\right)^{-1}(\0)$ 
is either empty or a codimension-$(n+1)$ submanifold of 
$S_1,$ and hence of $J^n\left(\mathfrak{X}^\infty(P)\right)$ as well.
We begin by computing $\Psi\left(j^{n}X(z)\right)$ in local coordinates for $P.$
Using superscripts to denote components of $X$ and subscripts to denote partial differentiation, using the summation convention for repeated indices,
and evaluating all derivatives at the same point $z=\gamma(0),$
we obtain:
\begin{align*}
\left(F\circ \gamma\right)'(0) &= F_i X^i \\
\left(F\circ \gamma\right)''(0) &= F_{ij} X^i X^j + F_i X^i_j X^j \\
&\vdots \\
\left(F\circ \gamma\right)^{k}(0) &= \textrm{(terms of lower order in $X$)} 
+ F_i X^i_{j_1j_2\cdots j_{k-1}} X^{j_1} X^{j_2}\cdots X^{j_{k-1}}
\end{align*}
We write the components of $\Psi$ as $\Psi_k,$ for $k=1,\dots,n+1,$ so
$\Psi_k(j^nX(z)) = \left(F\circ \gamma\right)^{(k)}(0).$
We now fix a $j^{n+1}X\left(z\right)\in S_1$ 
and consider $D\Psi$ at this point.
By definition of $S_1,$ we see that $z$ is neither a critical point of $X$ nor of $F.$
Since $X(z)\ne \0,$ there exists a $j\le n$ such that $X^j(z)\ne 0.$
Since $DF(z)\ne \0,$
there exists an $i\le n$ such that $F_i(z)\ne 0.$
It follows that, for every $k,$ 
the coefficient of $X^i_{jj\cdots j}$ in $\Psi_k$ is nonzero, where
$X^i_{jj\cdots j}$ means the $(k-1)^{\mathrm{th}}$-order partial derivative of $X^i$ with respect to 
the $j^{th}$ component of $z.$
This shows that the partial derivative of $\Psi_k$ with respect to $X^i_{jj\cdots j}(z)$ is nonzero.
Since $X^i_{jj\cdots j}$ does not appear in $\Psi_1,\cdots, \Psi_{(k-1)},$ and this 
argument holds for all $k=1,\dots,n+1$ and all $j^nX(z)\in S_1,$ this shows that
the restriction of $\Psi$ to $S_1$ is a submersion.
It follows that $S_2:=\left(\left.\Psi\right|_{S_1}\right)^{-1}(\0)$ is a smooth 
submanifold of 
$S_1,$ 
and hence of $J^n\left(\mathfrak{X}^\infty(P)\right)$ as well,
and that $S_2$ is either empty or of codimension $(n+1).$
Note that $\Psi^{-1}(\mathbf{0})\subseteq S_e\cup S_F \cup S_2.$
It follows that
\[
\left(\Psi\circ j^nX \right)^{-1}(\0) \subseteq \left(j^nX\right)^{-1}(S_e\cup S_F\cup S_2).
\]

Let $\mathcal{B}$ be the set of all $X\in \mathfrak{X}^\infty\left(P\right)$ such that
$j^nX$ is transverse to both $S_e$ and $S_2.$ 
By jet transversality for vector fields (Theorem \ref{jettransvf}), applied twice, 
$\mathcal{B}$ is a residual subset of $\mathfrak{X}^\infty\left(P\right).$
For any $X\in \mathcal{B},$ 
since $S_2$ is either empty or has codimension $n+1,$
its preimage $\left(j^nX\right)^{-1}(S_2)\subseteq P$ must be empty.
Similarly, since $S_e$ has codimension $n,$
its preimage $\left(j^nX\right)^{-1}(S_e)$ must be either empty or $0$-dimensional;
in other words it consists of isolated points (if any).
The set $\left(j^nX\right)^{-1}(S_F)$ consists of all $j^nX(z)$ with $z\in Z,$ where
$Z$ is a codimension-$1$ submanifold containing all of the critical points of $F.$
As noted earlier, any nontrivial solution of $X$ along which $F$ is constant
must remain in $\left(\Psi\circ j^nX \right)^{-1}(\0).$
Hence, the only such solutions (if any) are either equilibria,
which are isolated, or solutions remaining in $Z.$

The submanifold $Z$
may be covered with a countable number of
submanifold charts, in each of which 
$Z$ is a level set of some smooth function $f^i$
with no critical points.  
This function $f^i$ may be smoothly extended to a function on $P.$
Applying the above argument, with $f^i$ in place of $F,$
we conclude that there exists a residual subset $\mathcal{B}_i$ of  
$\mathfrak{X}^\infty\left(P\right)$ such that, for any $X\in\mathcal{B}_i,$ 
the only solutions (if any) of $X$ conserving $f^i$ are equilibria, and
these are isolated.
Define $\mathcal{C}= \mathcal{B}\cap \bigcap_i \mathcal{B}_i,$
and note that $\mathcal{C}$ is residual in $\mathfrak{X}^\infty\left(P\right).$
For any $X\in \mathcal{C},$
any solution of $X$ that remains in $Z$ must conserve
$f^i, $ for some $i$ and some time interval,
so it must be an equilibrium.
Hence the only solutions of $X$ along which $F$ is constant are equilibria,
and these are isolated.
\end{proof}

The previous theorem, with $F$ equal to the moment of inertia,
implies that, for generic vector fields,
there will be \emph{no} non-equilibrium constant-inertia solutions,
not even relative equilibria.
The existence of relative equilibria, 
for the Newtonian potential as for many others,
is of course related to the symmetry of the
of the vector field
 --- an $SO(2)$ symmetry in the case of the planar $N$-body problem.
Thus when generalising Saari's conjecture to arbitrary vector fields $X,$
it makes most sense to restrict attention to $SO(2)$-invariant vector fields,
for which the conjecture states:
``a solution of $X$ has constant moment of inertia
if and only if the solution is a relative equilibrium''.

As noted earlier, it suffices to consider only free symmetries
Let $G$ act freely, properly and smoothly on $P,$
so that $P/G$ is a smooth manifold.
If $F:P\to \R$ is $G$-invariant, then it descends to a function
$\bar{F}:P/G\to \R,$ the critical points of which are the
projections of the critical points of $F.$
If the critical points of $F$ are contained in some
$G$-invariant codimension-$1$ submanifold $Z,$
then it can be shown that $Z/G$
is a codimension-$1$ submanifold of $P/G.$
Applying Theorem \ref{gen1} to $P/G,$ we see that generic vector fields on
$P/G$ have no non-equilibrium solutions that conserve $\bar{F},$
and the equilibrium solutions are isolated.
The corresponding conclusion in the original phase space is:

\begin{theorem}\label{gen1sym}
Let $P$ be a
manifold and let
$G$ act freely and properly on $P.$
Let $\displaystyle{F:P\to\mathbf{R}}$ be a smooth $G$-invariant
 function such its critical points
are contained in a $G$-invariant codimension-$1$ submanifold of $P.$
Then there exists a residual subset $\mathcal{C}$ of $\mathfrak{X}^G(P)$
such that, for any vector field $X\in \mathcal{C},$
the only solutions of $X$ that conserve $F$ are relative equilibria,
and these are isolated (in the sense that they project to isolated points
in the reduced space).
\end{theorem}

The following corollary concerns the 
phase space $\R^{4N-4}$ of the $N$-body problem, and the
inertia function $I$ defined in the Introduction.
The corollary shows that 
a very ``direct'' generalisation of Saari's original conjecture
is generically true. 

\begin{corollary}
A solution to a generic 
$SO(2)$-symmetric vector field $X$ 
on $\R^{4N-4}$
has constant inertia if and only if the solution is a relative equilibrium.
\end{corollary}

\medskip

We now pursue analogous genericity results in smaller classes of vector fields: first Hamiltonian vector fields
on symplectic manifolds,
and then simple mechanical systems.
For any real valued function $H$ on a symplectic phase space $P$,
let $X_H$ be the associated Hamiltonian vector field on $P$.
The key ingredient is in the next theorem is
the fact that
the map $dH(z) \mapsto X_H(z)$, for any fixed $z,$ is surjective,
due to the nondegeneracy of the symplectic form.
We note that this is not true for general Hamiltonian vector fields on Poisson manifolds.

\begin{theorem}\label{Hamgen}
Let $\left(P,\omega\right)$ be a 
symplectic manifold
and let $\displaystyle{F:P\to\mathbf{R}}$ be a smooth
function such that the critical points of
$F$ are contained in a codimension-$1$ submanifold of $P.$
Then there exists a residual subset $\mathcal{C}$ of
$\mathcal{C}^\infty(P,\R)$
such that, for any $H\in \mathcal{C},$ the only solutions
to the Hamiltonian vector field $X_H$
along which $F$ is constant are equilibria, and these
are isolated.
\end{theorem}

\begin{proof}
Let $n$ be the dimension of $P.$
Define
\begin{align*}
\Phi: J^{n+1}(P,\R) &\longrightarrow  J^n\left(\mathfrak{X}^\infty(P)\right), \\
j^{n+1}H\left(z\right) &\longmapsto j^nX_H(z).
\end{align*}
The vector field $X_H$ is defined by $i_{X_H}\omega = dH$.
Since $\omega$ is smooth,
it follows that 
$\Phi$ is smooth.
For fixed $z,$ since $\omega(z)$ is nonsingular,
the map $dH(z)\mapsto X_H(z)$ is linear and surjective.
From this, it can be shown that $\Phi$ is surjective submersion.
Let $\Psi$ be defined as in the previous theorem, and recall that
any nontrivial solution of $X_H$ along which $F$ is constant
must remain in $\left(\Psi\circ j^nX_H \right)^{-1}(\0).$
Let $S_e,S_F,S_1$ and $S_2$ be as in the proof of the previous theorem,
and recall that $S_e$ is a submanifold of $J^n\left(\mathfrak{X}(P)\right)$ of codimension $n,$
$S_2$ is a submanifold of codimension $n+1,$
and $\Psi^{-1}(0)\subseteq S_e\cup S_F \cup S_2.$
Let $U_e=\Phi^{-1}\left(S_e\right)$ and
$U_2=\Phi^{-1}\left(S_2\right).$
Since $\Phi$ is a submersion, $U_e$
is a submanifold of $J^{n+1}(P,\R)$ of codimension $n$ and
$U_2$ is a submanifold of codimension $n+1.$

Let $\mathcal{B}$ be the set of all $H\in \mathcal{C}^\infty\left(P,\R\right)$ such that
$j^{n+1}H$ is transverse to both $U_e$ and $U_2.$
By jet transversality (Theorem \ref{jettrans}), applied twice,
$\mathcal{B}$ is a residual subset of $\mathcal{C}^\infty\left(P,\R\right).$
For any $H\in \mathcal{B},$
since $U_2$ is either empty or has codimension $n+1,$
its preimage $\left(j^{n+1}H\right)^{-1}(U_2)\subseteq P$ must be empty.
Similarly, since $U_e$ has codimension $n,$
its preimage $\left(j^{n+1}H\right)^{-1}(U_e)$ must be either empty or $0$-dimensional;
in other words it consists of isolated points (if any).
As noted earlier, any nontrivial solution of $X_H$ along which $F$ is constant
must remain in
$\left(\Psi\circ j^{n}X_H \right)^{-1}(\0)
=\left(\Psi\circ \Phi\circ j^{n+1}H \right)^{-1}(\0)
\subseteq \left(\Phi\circ j^{n+1}H \right)^{-1}\left(S_e\cup S_F \cup S_2\right).$
If we assume that $F$ has no critical points, then
this set is contained in \\
$\left(j^{n+1}H \right)^{-1}\left(U_e\cup U_2\right).$
Recalling that $\left(j^{n+1}H\right)^{-1}(U_2)$ is empty,
we conclude that
the only solutions along which $F$ is constant are equilibria, and these
are isolated.
It remains only to deal with the critical points of $F.$
This is done in the same manner as in the
previous theorem.
\end{proof}

We would like to prove an equivariant version of this result, analogous to
Theorem \ref{gen1sym}.
A natural approach is to use symplectic reduction (see \cite{AM78}),
and apply Theorem \ref{Hamgen}
on each reduced space.
There is a problem, however:
for each momentum value,
we have a different reduced space and a different reduced Hamiltonian.
One could show that, for every momentum value $\mu$, there
exists a residual set $\mathcal{C}_{\mu}$
of $G$-invariant Hamiltonians $H$ such that
$X_H$ has no constant-$F$ solutions of momentum $\mu.$
If we define $\mathcal{C}=\cap_\mu \mathcal{C}_\mu,$
then for any $H$ in $\mathcal{C},$ the vector field $X_H$ has no
constant-$F$ solutions; but since $\mu$ is a continuous variable,
$\mathcal{C}$ \emph{need not be residual.}
For this reason it seems reasonable to use Poisson reduction rather than symplectic reduction,
which would require an extension of Theorem \ref{Hamgen} to general Hamiltonian systems 
on Poisson manifolds.
Another possible approach is to avoid reduction entirely, and instead use
the theory of equivariant transversality (see \cite{F77}). 
This problems remains open.

\medskip

The next theorem concerns genericity in the class of Hamiltonian vector fields
with Hamiltonian of the form ``kinetic plus potential'',
for a fixed kinetic energy.
It is stated in the equivalent Lagrangian formulation.
For the proof of this theorem, standard jet transversality is insufficient,
since the potential is a function of configuration only.
For this reason we use the related 
version of transversality given in Theorem \ref{newtrans}.

\begin{theorem}\label{simpgen}
Let $Q$ be a 
Riemannian manifold and let
$\displaystyle{F:Q\to\mathbf{R}}$ be a smooth
 function such that the critical points of
$F$ are contained in codimension-$1$ submanifold of $Q.$
Let $K:TQ\to \R$ be the kinetic energy function defined by the given metric,
namely $K(q,\dot{q}) = \frac{1}{2}\left\| \dot q \right\|_q^2.$
Then there exists a residual subset $\mathcal{C}$ of $\mathcal{C}^\infty(Q,\R)$
such that for any $V\in \mathcal{C},$ the Euler-Lagrange equations for
the Lagrangian $L=K-V$
have no non-equilibrium solution along which $F$ is constant,
and the equilibrium solutions are isolated.
\end{theorem}

\begin{proof}
Let $n$ be the dimension of $Q$
and let $\pi:T^*Q\to Q$ be the cotangent bundle projection.
We consider Hamiltonians of the form $H=K+V.$ To be more precise,
for any function $V\in \mathcal{C}^\infty\left(Q,\R\right),$
let $\mathbb{F}L:TQ\to T^*Q$
be the Legendre transform of $L=K-V$ (which in fact depends only on the metric),
and define $H_V\in \mathcal{C}^\infty\left(T^*Q,\R\right)$
by $H_V=K\circ \left(\mathbb{F}L\right)^{-1} + V\circ \pi.$
Note that the solutions of the Euler-Lagrange equations for $L=K-V$
are the projections by $\pi$ of the solutions to the Hamiltonian
vector field $X_{H_V}$ on $T^*Q.$

Define 
$\Psi:\left(J^{2n-1}\left(\mathfrak{X}^\infty\left(P\right)\right)\right) \to \R^{2n+1}$ 
as in the proofs of the previous two theorems, namely 
$\Psi_k\left(j^{2n-1}X(z)\right)=\left(F\circ \gamma\right)^{(k)}(0),$
for every $X\in \mathfrak{X}^\infty(P)$ and every solution $\gamma$ of $X$ such that
$\gamma(0)=z.$
(Only $2n-1$ derivatives of $X$ are required to define $\Psi$  
because $F$ is a function of $Q$ only and $X$ is second-order.)
If $F\circ \gamma$ is constant,
then $\gamma(t)$ remains in $\left(\Psi \circ j^{2n-1}X\right)^{-1}(\0).$
Our strategy is to apply Theorem \ref{newtrans} with $\rho(V)=X_{H_V}.$
We thus define
\begin{align*}
\Phi: J^{2n}(Q,\R) \oplus P &\longrightarrow  J^{2n-1}\left(\mathfrak{X}^\infty(P)\right) \\
\left(j^{2n}V\left(q\right),z\right) &\longmapsto j^{2n-1}X_{H_V}(z)
\end{align*}
As noted in the proof of the previous theorem, 
$X_H(z)$ is a smooth function of $dH(z),$ for a fixed $z.$ 
Since $dH_V(z)$ is a smooth function of $dV(\pi(z)),$
it follows that $X_{H_V}(z)$ is a smooth function of $dV(\pi(z)),$ 
from which it follows that $\Phi$ is smooth.
We will show that $\Psi\circ \Phi$ is a submersion, except at points
corresponding to critical points of $F$ or $V.$

Let $U_e=\left\{ j^{2n}V(q) \in J^{2n}\left(Q,R\right)\, : \, dV\left(\pi(z)\right) = 0\right\}$
and $S_e = \Phi\left(U_e\oplus P\right).$
It is easily checked that $U_e$
is a closed codimension-$n$ subbundle of $J^{2n}\left(Q,R\right),$
which implies that $U_e\oplus P$ 
is a closed codimension-$n$ subbundle of $J^{2n}\left(Q,R\right)\oplus P.$
It is also easily checked that $\Phi$ is a proper injective immersion,
so its restriction to $U_e\oplus P$ is as well.
This implies that $S_e$ is a closed  codimension-$2n$ submanifold of
$J^{2n-1}\left(\mathfrak{X}^\infty(P)\right).$
Let $\mathcal{B}_e$ be the set of all $V\in \mathcal{C}^\infty(Q,\R)$ such that
$j^{2n}V$ is transverse to $U_e.$
By jet transversality (Theorem \ref{jettrans}), $\mathcal{B}_e$ is residual.

By assumption, the critical points of $F$ are contained in a codimension-$1$
submanifold $Z$ of $P.$
Let $S_F$ be the codimension-$1$ submanifold of $J^{2n-1}\left(\mathfrak{X}^\infty(P)\right)$
consisting of all $j^{2n-1}X(z)$ such that $z\in Z.$
Let $S_1$ be the complement of $S_e\cup S_F$ in $J^{2n-1}\left(\mathfrak{X}^\infty(P)\right),$
and note that $S_1$ is open dense.
Recall the map
$\Psi:\left(J^{2n-1}\left(\mathfrak{X}^\infty\left(\bar P
\right)\right)\right) \to \R^{2n+1}$ defined above. 
As in the proofs of the previous theorems, we can check that $\Psi$ is a submersion.
Let $S_2=\left(\left.G\right|_{S_1}\right)^{-1}(\0).$
Since $\Psi$ is a submersion,
$S_2$ is either empty or a codimension-$2n$
submanifold of $J^{2n-1}\left(\mathfrak{X}^\infty(P)\right).$
Note that $\left(\Psi\right)^{-1}(0)\subseteq S_e\cup S_F \cup S_2.$

We will show that $\left.\Phi\right|_{U_1}$ is transverse to $S_2.$
Since $\Psi$ is a submersion and $S_1$ is open dense in 
$J^n\left(\mathfrak{X}^\infty(P)\right),$ this is 
equivalent to showing that the
restriction of $\Psi\circ \Phi$ to $U_1:= \Phi^{-1}\left(S_1\right)$ 
is a submersion. 
We do so using canonical local coordinates 
$\left(z_1,\dots,z_{2n}\right)=\left(q^1,\dots,q^n,p_1,\dots,p_n\right)$
on $T^*Q.$
In these coordinates, 
$H=\frac{1}{2} p_i g^{ij} p_j + V,$
where $g^{ij}$ is the inverse of the metric tensor, and we use the
summation convention.
The vector field $X_{H_V}$ is defined by
$\left(X_{H_V}^1,\dots,X_{H_V}^{2n}\right) = 
\left(\dot{q}^1,\dots, \dot{q}^n, \dot{p}^1,\dots, \dot{p}^n \right)$ and 
\begin{align*}
\dot{q}^i &= \frac{\partial H}{\partial p_i} = g^{ij} p_j \\
\dot{p}_i &= - \frac{\partial H}{\partial q_i} = - \frac{1}{2} p_k g^{kl}_i p_l - V_i
\end{align*}
where subscripts denote differentiation with respect to variables $q^i.$
Given any path $\gamma(t)$ in $P,$ we can now compute
\begin{align*}
\left(F\circ \gamma\right)'(t) &= F_i \, \dot q^i = F_i \, g^{ij} \, p_j \\
\left(F\circ \gamma\right)''(t) 
&= \left(F_{il} \, g^{ij} + F_i \, g^{ij}_l \right) \dot q^l \, p_j+ F_i \, g^{ij} \, \dot p_j \\
&= \left(\textrm{terms with no $V_j$}\right) - F_i \, g^{ij} \, V_j \\
&\vdots \\
\left(F\circ \gamma\right)^{(k+1)}(t) 
&= \left(\textrm{terms of lower order in $V$}\right) 
+(-1)^k F_i \, g^{ij_1} \, V_{j_1j_2\dots j_k} 
\left(g^{j_2l}\, V_l\right) \dots \left(g^{j_kl}\, V_l\right)
\end{align*}
Note that $\left(\Psi\circ \Phi\right)_k= \left(F\circ \gamma\right)^{(k)}(t).$

We now fix a $\left(j^{2n}V(\pi(z)),z\right)\in U_1\oplus P$ 
and consider $D\left(\Psi\circ \Phi\right)$ at this point.
Since $\Phi \left(j^{2n}V(\pi(z)),z\right)$ is in $S_1,$ it is neither in $S_F$ nor
in $S_e.$
Thus we can assume that $\nabla F(z)$ and $\nabla V\left(\pi(z)\right)$ 
are both nonzero.
Thus there exists an $i$ and $j$ such that 
$F_i(z)$ and $V_j\left(\pi(z)\right)$ are both nonzero.
Since the metric tensor is positive definite, $g^{ii}(z)$ and $g^{jj}(z)$ are nonzero.
It follows that
$F_i \, g^{ii} \, V_{ij\dots j} \left(g^{jj}\, V_j\right)^{k-1},$
evaluated at our chosen $z,$ is nonzero, for
any $k\in \mathbf{N}.$
But this is the partial derivative of $\left(\Psi\circ \Phi\right)_k$ with respect to
$V_{ij\dots j}$ (the subscript $j$ is repeated $k-1$ times).
Since $V_{ij\dots j}$ doesn't appear in $\Psi_1\dots \Psi_{k-1},$
and this argument holds at any point in $U_1,$
we have shown that the restriction of $\Psi\circ \Phi$ to $U_1$ is a submersion.
It follows that $\left.\Phi\right|_{U_1}$ is transverse to $S_2.$
Since $S_2\subseteq S_1 = \Phi\left(U_1\right),$ this implies that
$\Phi$ is transverse to $S_2.$
Let $\mathcal{B}_2$ be the set of all $V\in \mathcal{C}^\infty\left(Q,\R\right)$
such that 
$X_{H_V}$ is transverse to $S_2.$
By Theorem \ref{newtrans}, with $\rho(V)=X_{H_V},$
$\mathcal{B}_2$ is a residual subset of $\mathcal{C}^\infty\left(Q,\R\right).$

Let $\mathcal{B}=\mathcal{B}_2\cup \mathcal{B}_e,$ which is also residual.
For any $V\in \mathcal{B},$
$j^{2n-1}X_{H_V}$ is transverse to $S_2$
and $j^{2n}V$ is transverse to $S_e.$
Since $j^{2n-1}X_{H_V}$ is transverse to $S_2,$
and $S_2$ is either empty or has codimension $2n,$
the preimage $\left(j^{2n-1}X_{H_V}\right)^{-1}\left(S_2\right)\subseteq T^*Q$
is empty.
Since and $j^{2n}V$ is transverse to $S_e,$ which has codimension $n,$
the set $\left(j^{2n} V\right)^{-1}\left(S_e\right)\subseteq Q$ has dimension $0,$
in other words consists of isolated points.
Recall that any solution to $X_{H_V}$ 
along which $F$ is constant must be entirely contained in
$\left(\Psi \circ j^{2n-1}X_{H_V}\right)^{-1}(\0),$ and
$\left(\Psi\right)^{-1}(0)\subseteq S_e\cup S_F \cup S_2.$
We have shown that $\left( j^{2n-1}X_{H_V}\right)\left(S_2\right)$ 
is empty. If we assume $F$ has no critical points, then 
any constant-$F$ solution must remain in 
$\left(j^{2n-1}X_{H_V}\right)^{-1}(S_e),$
which equals $\pi^{-1}\left(U_e\right).$
Since the elements of $U_e$ are isolated points of $Q,$ 
and our vector field is a second order equation on $Q,$
this means that any constant-$F$ solution is an equilibrium.

It remains only to deal with the critical points of $F.$ 
This is done in the same manner as in the 
previous two theorems.
\end{proof}

\begin{corollary} \label{saarigen}
For generic $N$-body potentials, there are no non-equilibrium constant-inertia trajectories.
\end{corollary}

We would like to prove an equivariant version of this theorem, analogous to Theorem \ref{gen1sym},
a corollary of which would be the following
generalisation of Saari's conjecture:
``a solution of the planar $N$-body problem,
 for a generic $SO(2)$-symmetric potential,
has constant inertia if and only if it is a relative equilibrium''.
As noted earlier, it suffices to consider 
free actions. Given a free action of $G$ 
on $Q,$ a $G$-invariant Riemannian metric on $Q$ (defining the
kinetic energy), and a $G$-invariant $F:Q\to \R,$ 
we would like to show that generic $G$-invariant potentials
have no constant-inertia solutions other than the relative equilibria.
A natural strategy is to try to use cotangent bundle reduction.
If we restrict attention to abelian groups $G$ then
the symplectic reduced spaces at different momentum values are all isomorphic to
$T^*\left(Q/G\right).$
The difficulty is that the reduced vector fields are momentum-dependent (see \cite{AM78}).
We have attempted to deal with this, using a generalisation of Theorem \ref{newtrans},
but so far without success.
Another (related) approach would be to use the symmetry-adapted coordinates given by the Cotangent bundle slice theorem, leading to the versions of the 
``bundle equations'' (``reconstruction equations'')
given in \cite{RSS05}.

\medskip
There are smaller classes of potentials relevant to Saari's conjecture
that we haven't yet addressed. One such class is the Newtonian potential for
arbitrary combinations of masses. 
Another is the class of potentials of the form
\[
V\left(\mathbf{q}_1,\dots,\mathbf{q}_N\right) 
= -\sum\limits_{1\le i<j\le N} m_i m_j f\left(\left\| \mathbf{q}_i - \mathbf{q}_j\right\|\right)
\]
with
\[
f(r)=\sum\limits_{1\le k \le 2N-1} \beta_k r^{\alpha_k}\, ,
\]
i.e., linear combinations of powers of
mutual distance.
(For example, this class includes the Lennard-Jones
potential.) 
While we see no theoretical obstacle to the use of our methods
to prove a genericity result in either of these classes,
we attest that the computational difficulty is considerable.

\begin{remark}
In the Newtonian $N$-body problem,
the relative equilibria correspond to central configurations.
The central configurations come in families scaled by inertia,
and the corresponding relative equilibria have angular momentum 
inversely proportional to the inertia (see \cite{Mey99}).
A well-known open problem asks: what combinations of masses admit
only a finite number of
central configurations, where ``finiteness'' is to be interpreted modulo
a scaling by inertia.
While our work does not address the question of finiteness,
we have shown that, for generic symmetric vector fields,
the relative equilibria are ``isolated'' (as points in the reduced space)
(Theorem \ref{gen1sym}).
However, since the relationship between relative equilibria and 
central configurations is specific to the Newtonian potential,
we will not be able to draw any conclusions about central configurations
unless and until a ``genericity of Saari's conjecture''
result is proven in the class of Newtonian
potentials with arbitrary masses.
\end{remark}

\section{Conclusion}
We have generalised Saari's conjecture in various ways, and shown that these generalised conjectures
are generically true. 
In Theorem \ref{genFsym} we showed that,
for any given $G$-invariant vector field $X$ (with $G$ acting freely), 
and for generic $G$-invariant functions $F$,
the only solutions to $X$ that conserve $F$ are the relative equilibria.
In Section \ref{sectgenX} we reversed the perspective, fixing $F$ and allowing $X$ to vary.
We have shown in Theorem \ref{gen1} that, 
for any given $G$-invariant $F$ without ``too many'' critical points,
and for generic $G$-invariant vector fields $X$, the only solutions to $X$
that conserve $F$ are the relative equilibria.
We then prove similar results for generic  Hamiltonian vector fields
(Theorem \ref{Hamgen}); and generic  Hamiltonian vector fields for Hamiltonians of the form ``kinetic plus potential'' with a given kinetic energy (Theorem \ref{simpgen}). 
However the latter two results do not address symmetry (equivalently,
we assume $G$ is trivial).
The problem of finding equivariant versions of Theorems \ref{Hamgen} and \ref{simpgen} remains open.
Our main tool in proving these results was jet transversality, 
including our new version in Theorem \ref{newtrans}.

Our genericity results seem ``natural'' and unsurprising from the point of view of transversality theory. 
Roughly speaking, by requiring $F\circ \gamma$ to be constant, for every solution $\gamma$ of $X,$
we are putting an infinite number of constraints on the derivatives of $F$ (if $X$ is fixed) or on
$X$ (if $F$ is fixed), and generic functions or vector fields will not satisfy these constraints.
The proofs of these theorems can be intricate, but the central idea  is consistent with a large body of 
``general position'' results.
Nonetheless, the consequences for Saari's conjecture are significant:
the lack of constant-inertia solutions 
(other than the relative equilibria) is not a property specific to the Newtonian potential,
or to the special relationship between moment of inertia and kinetic energy.
Instead, the non-conservation of
a function is a generic property in several classes of problems containing the specific one
for which Saari's original conjecture is formulated.
 
Of course, our result does not prove Saari's original conjecture; given the significance
of the $N$-body problem, this specific question is still of interest. However, our work does change the nature of future approaches to this conjecture.
In light of our genericity results, we would {\it expect} that the relative equilibria are the only constant-inertia solutions. If this is not the case for some potentials, then it is because those potentials have a non-generic relationship to the moment of inertia
(see, for example, \cite{San04}).
 In our experience, genericity is easiest to prove
in the most general classes of vector fields or functions. 
The more specific the class is, the more interesting
the result but the harder the proof: for example, Theorem \ref{simpgen} about generic potentials 
is more relevant to the $N$-body problem
than Theorem \ref{gen1} about generic vector fields, but it is harder to prove.
We expect that genericity within smaller classes, for example, the Newtonian potential with
generic masses, will be much harder,
for computational rather than theoretical reasons.
We speculate that any general proof
of Saari's original conjecture 
will be even more computationally difficult. Thus we no longer await the discovery
of an elegant conceptual proof of Saari's conjecture.

\bigskip
\textbf{Acknowledgments}
We thank Ian Melbourne and Francesco Fasso for their helpful suggestions.

\end{document}